\documentclass[12pt]{article}

\usepackage{amssymb,amsthm,amscd,array}

\let\ssection=\section
\renewcommand{\section}{\setcounter{equation}{0}\ssection}

\newcommand{\bbR}{\mathbb{R}}
\newcommand{\bbRP}{\mathbb{RP}}

\begin{document}



\def\a{\alpha}
\def\b{\beta}
\def\d{\delta}
\def\g{\gamma}
\def\om{\omega}
\def\r{\rho}
\def\s{\sigma}
\def\vfi{\varphi}
\def\vr{\varrho}
\def\l{\lambda}
\def\m{\mu}
\def\implies{\Rightarrow}

\oddsidemargin .1truein
\newtheorem{thm}{Theorem}[section]
\newtheorem{lem}[thm]{Lemma}
\newtheorem{cor}[thm]{Corollary}
\newtheorem{con}[thm]{Conjecture}
\newtheorem{pro}[thm]{Proposition}
\newtheorem{ex}[thm]{Example}
\newtheorem{rmk}[thm]{Remark}
\newtheorem{defi}[thm]{Definition}


\title{Projective geometry of polygons and discrete 4-vertex and 6-vertex theorems}

\author{
V.~Ovsienko\thanks{
CNRS, Centre de Physique Th\'eorique,
Luminy Case 907,
F--13288 Marseille, Cedex~9, FRANCE;
mailto:ovsienko@cpt.univ-mrs.fr
}
\and
S.~Tabachnikov\thanks{
Department of Mathematics,
University of Arkansas,
Fayetteville, AR 72701, USA;
mailto:serge@comp.uark.edu
}
}

\date{}

\maketitle

\bigskip

\thispagestyle{empty}

\begin{abstract}
The paper concerns discrete versions of the three well-known results of projective
differential geometry: the four vertex theorem,  the six affine vertex theorem and
the Ghys theorem on four zeroes of the Schwarzian derivative. We study geometry of closed
polygonal lines in $\bbRP^d$ and prove that polygons satisfying a certain convexity
condition have at least $d+1$ flattenings. This result provides a new approach to the above
mentioned classical theorems.
\end{abstract}

\newpage

\section{Introduction}

A vertex of a smooth plane curve is a point of its 4-th order contact with a circle (at a
generic point the osculating circle has 3-rd order contact with the curve). An affine vertex
(or sextactic point) of a smooth plane curve is a point of 6-th order contact with a conic. In
1909 S. Mukhopadhyaya
\cite{M} published the following celebrated theorems: {\it  every closed smooth convex
plane curve has at least 4 distinct vertices and at least 6 distinct affine vertices}. These
results generated a very substantial literature; from the modern point of view they are
related, among other subjects, to global singularity theory of wave fronts and Sturm theory
-- see e.g.
\cite{A1,A3,Bl,GMO,TU,U} and references therein.

A recent and unexpected result along these lines is the following theorem by E.~Ghys
\cite{G}: {\it the Schwarzian derivative of every projective line diffeomorphism has at
least 4 distinct zeroes} (see also \cite{OT,T1,DO}). The Schwarzian derivative vanishes
when the 3-rd jet of the diffeomorphism coincides with that of a projective transformation
(at a generic point a diffeomorphism can be approximated by a projective transformation up
to the 2-nd derivative). Ghys' theorem can be interpreted as the 4 vertex theorem in
Lorentzian geometry (cf. references above).

The goal of this note is to study polygonal analogs of the above three results.
In our opinion, such a discretization of smooth formulations is interesting for the
following reasons. First, a discrete theorem is a-priori stronger; it becomes, in the limit, a
smooth one, thus providing a new proof of the latter. An important feature of the discrete
approach is the availability of mathematical induction which can considerably simplify
the proofs. Second, the very
operation of discretization is non-trivial: a single smooth theorem may lead to
non-equivalent discrete ones. An example of this phenomenon is provided by two recent
versions of the 4 vertex theorem for convex plane polygons \cite{S,S1,W,T2} -- see Remark
\ref{FourRmk} below. To the best of our knowledge, these results are the only available
discrete versions of the 4 vertex theorem.

In this regard we would like to attract attention to the celebrated Cauchy lemma~(1813): 
{\it given two convex plane (or spherical) polygons whose respective sides are congruent,
the cyclic sequence of the differences of the respective angles of the polygons changes
sign at least 4 times}. This result plays a crucial role in the proof of convex polyhedra
rigidity (see
\cite{C} for a survey). The Cauchy lemma implies, in the limit, the smooth 4 vertex
theorem and can be considered as the first result in the area under discussion.

\goodbreak

\section{Theorems on plane polygons}\label{PlaneSection}

In this section we formulate our results for plane polygonal curves. The proofs will be
given in Section \ref{ProofSection}.

\subsection{Discrete 4 vertex theorem}\label{FourSection}

The osculating circle of a smooth plane curve is a circle that has 3-rd order contact
with the curve. One may say that the osculating circle passes through 3 infinitely close
points; at a vertex the osculating circle passes through 4 infinitely close points.
Moreover, a generic curve crosses the osculating circle at a generic point and stays on one
side of the osculating circle at a vertex. This well-known fact motivates the following
definition.

Let $P$ be a plane convex $n$-gon; we assume that $n\geq4$ throughout this section. Denote
the consecutive vertices by
$V_1,\ldots,V_n$, where we understand the indices cyclically, that is, $V_{n+1}=V_1$, etc.

\begin{defi}
{\rm
A triple of vertices $(V_i,V_{i+1},V_{i+2})$ is called extremal if $V_{i-1}$ and $V_{i+3}$
lie on the same side of the circle through $V_i,V_{i+1},V_{i+2}$ (this does not exclude the
case when $V_{i-1}$ or $V_{i+3}$ belongs to the circle)\footnote{We have a terminological
difficulty here: dealing with polygons, we cannot use the term ``vertex'' in the same sense
as in the smooth case; thus the term ``extremal''.}.
}
\label{VertexDef}
\end{defi} 

\null

\vspace{2cm}

~\hspace{0,1cm}
\special{illustration 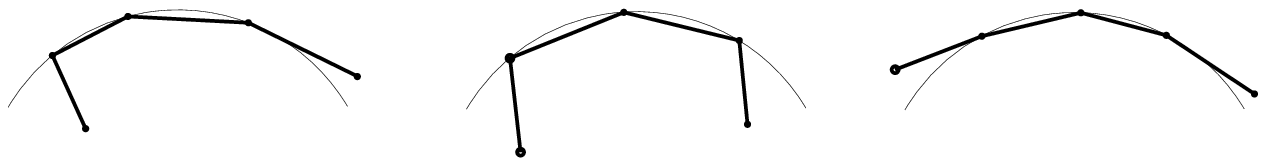}

\vspace{0,2cm}
\hskip 0,7cm
a) not extremal
\hskip 4,5cm
b) extremal

\bigskip

\centerline{\sc Figure 1}

\bigskip

The next result follows from a somewhat more general theorem due to O.~Musin and
V.~Sedykh \cite{S} (see also \cite{S1}).

\goodbreak

\begin{thm}
\label{FourTheorem}
Every plane convex polygon $P$ has at least 4 extremal triples of vertices.
\end{thm}

\begin{ex}
\label{FourEx}
{\rm
If $P$ is a quadrilateral then the theorem holds tautologically since $(i-1)$-st
vertex coincides with $(i+3)$-rd for every $i$.
}
\end{ex}
 
\begin{rmk}
\label{FourRmk}
{\rm
An alternative approach to discretization of the 4 vertex theorem consists in
inscribing circles in consecutive triples of sides of a polygon (the center of such a circle
is the intersection point of the bisectors of the consecutive angles of the polygon). Then a
triple of sides $(\ell_i,\ell_{i+1},\ell_{i+2})$, is called extremal if the lines
$\ell_{i-1},\ell_{i+3}$ either both intersect the corresponding circle or both do not
intersect it. With this definition an analog of Theorem \ref{FourTheorem} holds \cite{W,T2},
 also providing, in the limit, the smooth 4 vertex theorem.

Both formulations, concerning circumscribed and inscribed circles, make sense on the
sphere and, moreover, are equivalent via projective duality.
}
\end{rmk}

\subsection{Discrete 6 affine vertex theorem}\label{SixSection}

Five points in the plane determine a conic. Considering the plane as an affine part of the
projective plane, the complement of the conic has two connected components. Let
$P$ be a plane convex
$n$-gon; we assume that $n\geq6$ throughout this section. Similarly to the preceding
section, we give the following definition.

\begin{defi}
{\rm
Five consecutive vertices $V_i,\ldots,V_{i+4}$ are called extremal if $V_{i-1}$ and
$V_{i+5}$ lie on the same side of the conic through these 5 points (this does not
exclude the case when $V_{i-1}$ or $V_{i+5}$ belongs to the conic).
}
\label{AffVertexDef}
\end{defi} 

If $P$ is replaced by a smooth convex curve and $V_i,\ldots,V_{i+4}$ are infinitely close
points we recover the definition of the affine vertex. The following theorem is, therefore,
a discrete version of the smooth 6 affine vertex theorem.

\begin{thm}
\label{SixTheorem}
Every plane convex polygon $P$ has at least 6 extremal quintuples of vertices.
\end{thm}

\begin{ex}
\label{SixEx}
{\rm
If $P$ is a hexagon then the theorem holds tautologically for the same reason as in Example
\ref{FourEx}. 
}
\end{ex}

\begin{rmk}
\label{SixRmk}
{\rm
Interchanging sides and vertices and replacing  circumscribed conics by inscribed ones, we
arrive at a ``dual'' theorem. The latter is equivalent to Theorem \ref{SixTheorem} via
projective duality -- cf. Remark \ref{FourRmk}. 
}
\end{rmk}

\subsection{Discrete Ghys theorem}\label{GhysSection}

A discrete object of study in this section is a pair of cyclically ordered $n$-tuples
$X=(x_1,\ldots,x_n)$ and
$Y=(y_1,\ldots,y_n)$ in $\bbRP^1$ with $n\geq4$. Choosing an orientation of $\bbRP^1$, we
assume that the cyclic order of each of the two $n$-tuples is induced by the orientation.

Recall that an ordered quadruple of points in $\bbRP^1$ determines a number,
called the cross-ratio, which is projectively invariant. Choosing an affine parameter so
that the points are given by real numbers $a<b<c<d$, the cross-ratio is
\begin{equation}
[a,b,c,d]=
\frac{(c-a)(d-b)}{(b-a)(d-c)}.
\label{CrossEq}
\end{equation}

\begin{defi}
{\rm
A triple of consecutive indices $(i,i+1,i+2)$ is called extremal if
the difference of cross-ratios
\begin{equation}
[y_j,y_{j+1},y_{j+2},y_{j+3}]-[x_{j},x_{j+1},x_{j+2},x_{j+3}]
\label{CrossVertex}
\end{equation}
changes sign as $j$ changes from $i-1$ to $i$ (this does not exclude the case when either
of the differences vanishes).
}
\label{GhysDef}
\end{defi}

\begin{thm}
For every pair of $n$-tuples of points $X,Y$ as above there exist at least four extremal
triples.
\label{GhysTheorem}
\end{thm}

\begin{ex}
{\rm
If $n=4$ then the theorem holds for the following simple reason. The cyclic permutation
of four points induces the next transformation of the cross-ratio
\begin{equation}
[x_4,x_1,x_2,x_3]=
\frac{[x_1,x_2,x_3,x_4]}{[x_1,x_2,x_3,x_4]-1}
 \label{Transform}
\end{equation}
which is an involution.
Furthermore, if $a>b>1$ then $a/(a-1)<b/(b-1)$. Therefore, each triple of indices is
extremal.
}
\end{ex}

Let us interpret Theorem \ref{GhysTheorem} in geometrical terms similarly to Theorems
\ref{FourTheorem} and \ref{SixTheorem}. There exists a unique
projective transformation that takes $x_i,x_{i+1},x_{i+2}$ to $y_i,y_{i+1},y_{i+2}$,
respectively. Consider the graph $G$ of this transformation as a curve in
$\bbRP^1\times\bbRP^1$; the three points
$(x_i,y_i),(x_{i+1},y_{i+1}),(x_{i+2},y_{i+2})$ lie on this graph.
An ordered couple of points $(x_j,x_{j+1})$ in oriented $\bbRP^1$ defines the unique
segment. An ordered couple of points $((x_j,y_j),(x_{j+1},y_{j+1}))$ in
$\bbRP^1\times\bbRP^1$ also defines the unique segment, the one whose projection on
each factor is the defined segment in $\bbRP^1$. 
The triple $(i,i+1,i+2)$ is extremal if and
only if the topological index of intersection of the broken line 
$(x_{i-1},y_{i-1}),\ldots,(x_{i+3},y_{i+3})$ 
with the graph $G$ is zero. This fact
can be checked by a direct computation using (\ref{CrossEq}) that we omit.

\null

\vspace{4cm}

~\hspace{3cm}\special{illustration 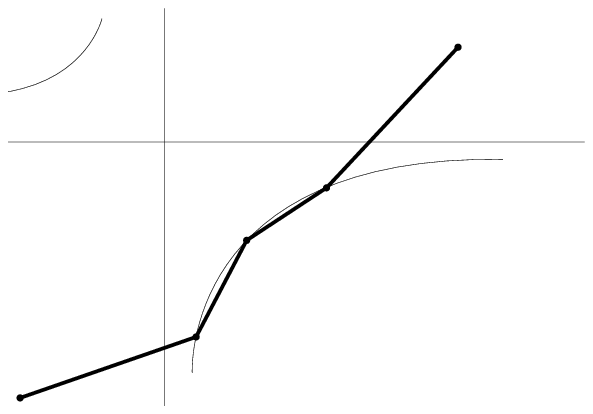}

\vspace{0,2cm}

\bigskip

\centerline{\sc Figure 2}

\bigskip

Let us also comment on the relation between Definition \ref{GhysDef} and zeroes of the
Schwar\-zian derivative of a diffeomorphism of the projective line. Let
$$
x_0=0,\;x_1=\varepsilon,\;x_2=2\varepsilon,\;x_3=3\varepsilon
$$
be four infinitely close
points given in an affine coordinate, and let $y_i=f(x_i)$ where~$f$ is a diffeomorphism of
$\bbRP^1$. Then
$$
[y_0,y_1,y_2,y_3]-[x_0,x_1,x_2,x_3]=
\varepsilon^2S(f)(0)+O(\varepsilon^3)
$$
where 
$$
S(f)=\frac{f'''}{f'}-
\frac{3}{2}\left(\frac{f''}{f'}\right)^2
$$
is the Schwarzian derivative of $f$. Thus, for $\varepsilon\to0$, Definition \ref{GhysDef}
corresponds to vanishing of the Schwarzian derivative.

\section{Main Theorem}\label{Main}

All the theorems from Section \ref{PlaneSection} are consequences of one theorem on the
least number of flattenings of a closed polygon in real projective space. 

In his remarkable work \cite{Ba}, M.~Barner introduced the notion of a {\it strictly convex}
curve in real projective space: this is a smooth closed curve $\g\subset\bbRP^d$ such that
for every $(d-1)$-tuple of points of $\g$ there exists a hyperplane through these points
that does not intersect $\g$ at other points. Barner discovered the following theorem: {\it a
strictly convex curve has at least $d+1$ distinct flattening points.} Recall that a flattening
point of a projective space curve is a point at which the osculating hyperplane is
stationary; in other words, this is a singularity of a projectively dual curve. In fact,
Barner's result is considerably stronger but we will not dwell on it here -- see \cite{T1}
for an exposition.

Our goal in this section is to provide a discrete version of Barner's theorem. First we need
to develop an elementary intersection formalism for polygonal lines.

\subsection{Multiplicities of intersections}\label{Formalism}

Throughout this section we will be considering closed polygons $P\subset\bbRP^d$ with
vertices $V_1,\ldots,V_n$, $n\geq{}d+1$, in general position. This means that
for every set of vertices
$V_{i_1},\ldots,V_{i_k}$ where
$k\leq{}d+1$ the span of $V_{i_1},\ldots,V_{i_k}$ is $(k-1)$-dimensional.

\begin{defi}
\label{TransDefi}
{\rm
A polygon $P$ is said to be transverse to a hyperplane $H$ at point $X\in{}P\cap{}H$ if
(a) $X$ is an interior point of an edge and this edge is transverse to~$H$, or (b) $X$ is a
vertex, the two edges incident to $X$ are transverse to $H$ and are locally separated by
$H$.
}
\end{defi}
Clearly, transversality is an open condition.
\begin{defi}
{\rm
A polygon $P$ is said to intersect a hyperplane $H$ with multiplicity~$k$ if
for every hyperplane $H'$ sufficiently close to $H$ and transverse to $P$, the number of
points $P\cap{}H'$ does not exceed $k$ and, moreover, $k$ is achieved for some $H'$.
}
\label{MultDefi}
\end{defi}

\null

\vspace{2cm}

~\hspace{1,5cm}\special{illustration 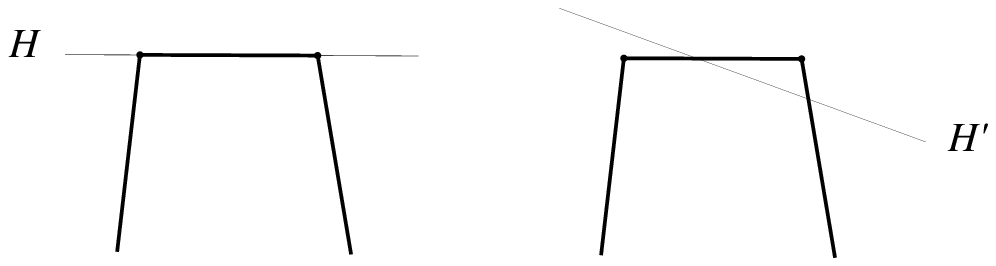}

\vspace{0,2cm}

\bigskip

\centerline{\sc Figure 3}

\bigskip

This definition does not exclude the case when a number of vertices of $P$ lie in~$H$.

\begin{lem}
\label{MultLem}
Let $V_{i_1},\ldots,V_{i_k}$ with $k\leq{}d$ be vertices of $P$. Then any hyperplane $H$
passing through $V_{i_1},\ldots,V_{i_k}$ intersects $P$ with multiplicity at least $k$.
\end{lem}
\begin{proof}
Move each $V_{i_j}$, $j=1,\ldots,k$ slightly along the
edge $(V_{i_j},V_{i_j+1})$ to obtain a new point $V_{i_j}'$. Let us
show that a generic hyperplane $H'$ through $V'_{i_1},\ldots,V'_{i_k}$ is transverse to
$P$. This will imply the lemma because $H'$ has at least $k$ intersections with $P$.

It suffices to show that $H'$ contains no vertices of $P$. 
Note first that since $P$ is in general
position, a generic hyperplane
$H$ through $V_{i_1},\ldots,V_{i_k}$ does not contain any other vertex. 
The same holds true for
every hyperplane sufficiently close to~$H$. It remains to show that the chosen $H'$ does not
contain either of $V_{i_1},\ldots,V_{i_k}$. 

Let $H'$ contain $V_{i_j}$; then $H'$ also contains the edge $(V_{i_j},V_{i_j+1})$ and
therefore~$V_{i_j+1}$. If
$i_j+1\not\in\{i_1,\ldots,i_k\}$ we obtain a contradiction with the previous paragraph. On
the other hand, if $i_j+1\in\{i_1,\ldots,i_k\}$ we can continue in the same way. However,
we cannot continue indefinitely since $k<n$.
\end{proof}

The next definition is topological in nature.
\begin{defi}
{\rm
Consider a continuous curve in $\bbRP^d$ with endpoints $A$ and $Z$ and let~$H$ be a
hyperplane not containing $A$ and $Z$. We say that $A$ and $Z$ are on one side of $H$ if one
can connect $A$ and $Z$ by a curve not intersecting $H$ such that the obtained closed curve
is contractible; and $A$ and $Z$ are separated by $H$ otherwise.
 }
\label{SeparateDefi}
\end{defi}
Clearly, if one only has two points $A$ and $Z$ (and no curve connecting them), then one
cannot say that the points are on the one side of or separated by a hyperplane.

\begin{lem}
\label{SeparateLem}
Let $\Gamma=(A,\ldots,Z)$ be a broken line in $\bbRP^d$ in general position and let
$H$ be a hyperplane not containing $A$ and $Z$. Denote by $k$ the multiplicity
of the intersection of
$\Gamma$ with $H$. Then $A$ and $Z$ are separated by $H$ if $k$
is odd and not separated otherwise.
\end{lem}
\begin{proof}
Connect $Z$ and $A$ by a segment to obtain a closed polygon $\overline{\Gamma}$ and
consider a hyperplane $H'$ close to $H$, transverse to
$\overline{\Gamma}$ and intersecting $\Gamma$ in $k$ points. Since
$\overline{\Gamma}$ is contractible, $H'$ intersects
$\overline{\Gamma}$ in an even number of points. Therefore, $H'$ intersects
the segment $(Z,A)$ for odd $k$ and does not intersect it for even $k$.
\end{proof}

\goodbreak

The next definition introduces a significant class of polygons which is our main object of
study.
\begin{defi}
{\rm
A polygon $P$ is called strictly convex if through every $d-1$ vertices there passes a
hyperplane~$H$ such that the multiplicity of its intersection with~$P$ equals $d-1$.
 }
\label{MainDefi}
\end{defi}
The preceding definition becomes, in the smooth limit, that of strict convexity of smooth
curves, due to Barner.

\goodbreak

\begin{defi}
\label{FlatDefi}
{\rm
A $d$-tuple of consecutive vertices $(V_i,\ldots,V_{i+d-1})$ of a polygon
$P$ in~$\bbRP^d$ is called a flattening if the endpoints $V_{i-1}$ and $V_{i+d}$ of the
broken line $(V_{i-1},\ldots,V_{i+d})$ are: (a) separated by the hyperplane through
$(V_i,\ldots,V_{i+d-1})$ if $d$ is even, (b) not separated if $d$ is odd.
 }
\end{defi}

\null

\vspace{2,5cm}

~\hspace{2cm}\special{illustration 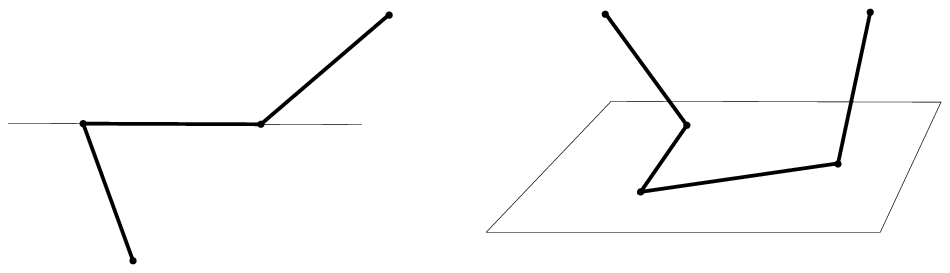}

\vspace{0,2cm}

\bigskip

\centerline{\sc Figure 4}

\bigskip

\begin{rmk}
\label{FlatRmk}
{\rm
A curve in $\bbRP^d$ can be lifted to $\bbR^{d+1}\setminus\{0\}$; this lifting is not
unique. Given a polygon $P\subset\bbRP^d$ with vertices $V_1,\ldots,V_n$, we lift it to
$\bbR^{d+1}$ as a polygon $\widetilde{P}$ and denote its vertices by
$\widetilde{V}_1,\ldots\widetilde{V}_n$.
Then a $d$-tuple $(V_i,\ldots,V_{i+d-1})$ is a flattening if and only if the determinant
\begin{equation}
\label{DeltaEq}
\Delta_j
=\left|\widetilde{V}_j,\ldots\widetilde{V}_{j+d}\right|
\end{equation}
changes sign as $j$ changes from $i-1$ to $i$.
}
\end{rmk}
Note that this property is independent of the lifting.

\subsection{Simplex is strictly convex}\label{SimplSect}

Define a simplex $S_d\subset\bbRP^d$ with vertices $V_1,\ldots,V_{d+1}$ as the
projection from the punctured~$\bbR^{d+1}$ of the polygonal line:
\begin{equation}
\label{SimplEq}
\widetilde{V}_1=(1,0,\ldots,0), \;
\widetilde{V}_2=(0,1,0,\ldots,0),\;
\ldots,
\widetilde{V}_{d+1}=(0,\ldots,0,1)
\end{equation}
and
\begin{equation}
\label{ParityEq}
\widetilde{V}_{d+2}=(-1)^{d+1}\widetilde{V}_1.
\end{equation}
The last vertex has the same projection as the first one;  $S_d$ is contractible for odd $d$
and non-contractible for even $d$.

\null

\vspace{2cm}

~\hspace{2cm}\special{illustration 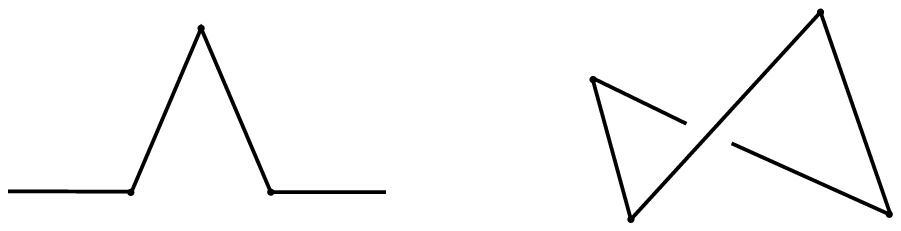}

\vspace{0,2cm}

\hskip 3cm
a) $d=2$
\hskip 5cm
b) $d=3$

\bigskip

\centerline{\sc Figure 5}

\bigskip

\begin{pro}
The polygon $S_d$ is strictly convex.
\label{SymplPro}
\end{pro}
\begin{proof}
We need to prove that through every $(d-1)$-tuple
$(V_1,\ldots,\widehat{V_i},\ldots,\widehat{V_j},\ldots,V_{d+1})$ there passes a
hyperplane
$H$ intersecting $P$ with multiplicity $d-1$. 
Choose a point $W$ on the line $(\widetilde{V}_i,\widetilde{V}_j)$ so that $W$ lies on the
segment $(\widetilde{V}_i,\widetilde{V}_j)$ if $j-i$ is even and does not lie on this
segment if $j-i$ is odd. Define $\widetilde{H}$ as the linear span of
$\widetilde{V}_1,\ldots,\widehat{\widetilde{V}_i},
\ldots,\widehat{\widetilde{V}_j},\ldots,\widetilde{V}_{d+1},W$. We claim that its
projection $H\subset\bbRP^d$ intersects $S_d$ with multiplicity $\leq{}d-1$.

Let $H'$ be a hyperplane close to $H$ and transverse to $S_d$; assume, moreover, that $H'$
contains no vertices. It suffices to show that $H'$ cannot intersect $S_d$ in more than
$d-1$ points. On the one hand, $H'$ cannot intersect all the edges of $S_d$. 
If it were the case, then $\widetilde{H'}$ would separate all pairs of consecutive vertices
which contradicts the choice of $W$. On the other hand, if the number of intersections of
$H'$ and $S_d$ were greater than $d-1$, it would be equal to $d+1$ since, for topological
reasons, the parity of this number of intersections is that of $d-1$. We obtain a
contradiction which proves the claim.

Finally, by Lemma \ref{MultLem}, the multiplicity of the intersection of $H$ with $S_d$ is
not less than $d-1$.
\end{proof}

\goodbreak

A curious property of the simplex is that each of its $d$-tuples of vertices is a flattening.

\begin{lem}
\label{BasaLem}
Simplex $S_d$ has $d+1$ flattenings.
\end{lem}
\begin{proof}
The determinant (\ref{DeltaEq}) involves all $d+1$ vectors
$\widetilde{V}_1,\ldots,\widetilde{V}_{d+1}$. If $d$ is odd then, according to
(\ref{ParityEq}),
$\widetilde{V}_{d+2}=\widetilde{V}_1$, and Definition \ref{FlatDefi} concerns the cyclic
permutation  of the vectors that changes sign of the determinant.
On the other hand, if $d$ is even then 
$\widetilde{V}_{d+2}=-\widetilde{V}_1$ which also leads to the change of sign in
(\ref{DeltaEq}).
\end{proof}

\subsection{Barner theorem for polygons}\label{BarnerSect}

Now we formulate the result which serves the main technical tool in the proof of Theorems
\ref{FourTheorem}, \ref{SixTheorem} and \ref{GhysTheorem}. Recall that we consider
generic polygons in
$\bbRP^d$ with at least $d+1$ vertices.

\begin{thm}
\label{BarnerThm}
A strictly convex polygon in $\bbRP^d$ has at least $d+1$ flattenings.
\end{thm}
\begin{proof}
Induction in the number of vertices $n$.

Induction starts with $n=d+1$. Up to projective transformations, the unique
strictly convex $(d+1)$-gon is the simplex $S_d$. Indeed, every generic $(d+1)$-tuple of
points in $\bbRP^d$ can be taken to any other by a projective transformation.
Therefore, all generic broken lines with $d$ edges are projectively equivalent. It remains to
connect the last point with the first one, and there are exactly two ways to do it. One of
them gives a contractible polygon and the other non-contractible one. One of these polygons
is~$S_d$, while the other cannot be strictly convex since the parity of its intersections
with a hyperplane is opposite to $d-1$. The base of induction is then provided by
Lemma~\ref{BasaLem}.

Let $P$ be a strictly convex $(n+1)$-gon with vertices $V_1,\ldots,V_{n+1}$. Delete
$V_{n+1}$ and connect $V_n$ with $V_1$ in such a way that the new edge $(V_n,V_1)$
together with the two deleted ones, $(V_n,V_{n+1})$ and $(V_{n+1},V_1)$, forms a
contractible triangle. Denote the new polygon by $P'$.

Let us show that $P'$ is strictly convex. $P$ is strictly convex, therefore through every
$d-1$ vertices of $P'$ there passes a hyperplane $H$ intersecting $P$ with multiplicity
$d-1$. We want to show that the multiplicity of the intersection of $H$ with $P'$ is also
$d-1$. Let $H'$ be a hyperplane close to $H$ and transverse to $P$ and $P'$. The
number of intersections of $H'$ with $P'$ does not exceed that with $P$. Indeed, if $H'$
intersects the new edge, then it intersects one of the deleted ones since the triangle in
contractible.

By the induction assumption, $P'$ has at least $d+1$ flattenings. To prove the theorem, it
remains to show that $P'$ cannot have more flattenings than $P$.

Consider the sequence of determinants (\ref{DeltaEq})
$\Delta_1,\Delta_2,\ldots,\Delta_{n+1}$. Replacing $P$ by~$P'$ we remove $d+1$
consecutive determinants
\begin{equation}
\Delta_{n-d+1},\Delta_{n-d+2},\ldots,\Delta_{n+1}
\label{DeltaOldEq}
\end{equation}
 and add in their stead $d$ new
determinants 
\begin{equation}
\Delta'_{n-d+1},\Delta'_{n-d+2},\ldots,\Delta'_n
\label{DeltaNewEq}
\end{equation}
 where
\begin{equation}
\Delta'_{n-d+i}=
\left|
\widetilde{V}_{n-d+i}\ldots\widehat{\widetilde{V}}_{n+1}\ldots\widetilde{V}_{n+i+1}
\right|
\label{DeltaPrimeEq}
\end{equation}
with $i=1,\ldots,d$. The transition from (\ref{DeltaOldEq}) to (\ref{DeltaNewEq}) is done
in two steps. First, we add (\ref{DeltaNewEq}) to (\ref{DeltaOldEq}) so that the two
sequences alternate, that is, we put $\Delta'_j$ between $\Delta_j$ and $\Delta_{j+1}$.
Second, we delete the ``old'' determinants (\ref{DeltaOldEq}). We will prove that the first
step preserves the number of sign changes while the second step obviously cannot increase
this number.

\begin{lem}
\label{SignLem}
If $\Delta_{n-d+i}$ and $\Delta_{n-d+i+1}$ are of the same sign, then $\Delta'_{n-d+i}$
is of the same sign too.
\end{lem}
\noindent
\textit{Proof of the lemma}. Since $P$ is in general position,
the removed vector $\widetilde{V}_{n+1}$ is a linear combination of $d+1$ vectors
$\widetilde{V}_{n-d+i},\ldots,\widetilde{V}_n,
\widetilde{V}_{n+2},\ldots,\widetilde{V}_{n+i+1}$:
\begin{equation}
\widetilde{V}_{n+1}=
a\widetilde{V}_{n-d+i}+
b\widetilde{V}_{n+i+1}+\cdots
\label{DeltaCombEq}
\end{equation}
where dots is a linear combination of the rest of the vectors.
It follows from (\ref{DeltaPrimeEq}) that 
\begin{equation}
\Delta_{n-d+i}=(-1)^{i-1}b\Delta'_{n-d+i},
\quad
\Delta_{n-d+i+1}=(-1)^{d-i}a\Delta'_{n-d+i}.
\label{ExpressEq}
\end{equation}

It is time to make use of strict convexity of $P$. Let $H$ be a hyperplane in $\bbRP^d$
through $d-1$ vertices $V_{n-d+i+1},\ldots,\widehat{V}_{n+1},\ldots,V_{n+i}$ that
intersects $P$ with multiplicity $d-1$, and let $\widetilde{H}$ be its lifting to
$\bbR^{d+1}$. Choose a linear function $\varphi$ in $\bbR^{d+1}$ vanishing on
$\widetilde{H}$ and such that $\varphi(\widetilde{V}_{n+1})>0$. We claim that
\begin{equation}
(-1)^{d-i}\varphi(\widetilde{V}_{n-d+i})>0
\quad\hbox{and}\quad
(-1)^{i-1}\varphi(\widetilde{V}_n)>0.
\label{DeltaExpressEq}
\end{equation}
Indeed, by Lemma \ref{MultLem}, the multiplicity of the intersection of $\widetilde{H}$
with the polygonal lines
$(\widetilde{V}_{n-d+i},\ldots,\widetilde{V}_{n+1})$ and
$(\widetilde{V}_{n+1},\ldots,\widetilde{V}_{n+i+1})$ are at least
$d-i$ and $i-1$, respectively. Since $H$ intersects $P$ with multiplicity $d-1$, the two
above multiplicities are indeed equal to $d-i$ and $i-1$.
The inequalities (\ref{DeltaExpressEq}) now readily follow from Lemma \ref{SeparateLem}.

Finally, evaluate $\varphi$ on (\ref{DeltaCombEq}):
$$
\varphi(\widetilde{V}_{n+1})=
a\varphi(\widetilde{V}_{n-d+i})+
b\varphi(\widetilde{V}_{n+i+1}).
$$
It follows from (\ref{DeltaExpressEq}) and $\varphi(\widetilde{V}_{n+1})>0$ that at least
one of the numbers $(-1)^{i-1}b$ or $(-1)^{d-i}a$ is positive. In view of (\ref{ExpressEq}),
Lemma \ref{SignLem} follows.
\end{proof}

Theorem \ref{BarnerThm} is proved.

\begin{rmk}
{\rm
Strict convexity is necessary for the existence of~${d+1}$ flattenings. One can
easily construct a closed polygon without any flattenings and even $C^0$-approximate an
arbitrary closed smooth curve by such polygons. The smooth case such approximation is
well known: given a curve $\g_0$, the approximating one, $\g$, spirals around in its tubular
neighborhood. In the polygonal case we take a sufficiently fine straightening of $\g$.
}
\end{rmk}

\section{Applications of the main theorem}\label{Last}

\subsection{Proof of Theorems \ref{FourTheorem}, \ref{SixTheorem}
and \ref{GhysTheorem}}\label{ProofSection}

Now we prove the results announced in Section \ref{PlaneSection}.
The idea is the same in all three cases and is precisely that of Barner's proof of
smooth versions of these theorems -- see~\cite{Ba} and also~\cite{T1}. We will consider in
detail Theorem
\ref{SixTheorem} indicating the necessary changes in the two other cases.

Let $P$ be as in Theorem \ref{SixTheorem}.
Consider the Veronese map ${\cal V}:\bbRP^2\to\bbRP^5$ of degree 2:
\begin{equation}
\label{VeroneseMap}
{\cal V}:(x:y:z)\mapsto
(x^2:y^2:z^2:xy:yz:zx).
\end{equation}
The image ${\cal V}(P)$ is a piecewise smooth curve. Homotop every edge to a straight
segment, keeping the endpoints ${\cal V}(V_i),{\cal V}(V_{i+1})$ fixed, to obtain a polygon
$Q$ in $\bbRP^5$. Assume first that $Q$ is in general position.
\begin{lem}
\label{PeregLem}
A quintuple ${\cal V}(V_i),\ldots,{\cal V}(V_{i+4})$ is a flattening of $Q$, if and only if
$(V_i,\ldots,V_{i+4})$ is an extremal quintuple of vertices of $P$.
\end{lem}
\begin{proof}
The Veronese map establishes a one-to-one correspondence between conics in $\bbRP^2$
and hyperplanes in $\bbRP^5$: the image of the conic is the intersection of a hyperplane
with the quadratic surface ${\cal V}(\bbRP^2)$. Since ${\cal V}$ is an
embedding, the points $V_{i-1}$ and $V_{i+5}$ lie on one side of the conic through
$(V_i,\ldots,V_{i+4})$ if and only if the points ${\cal V}(V_{i-1})$ and ${\cal V}(V_{i+5})$
lie on one side of the corresponding hyperplane.
\end{proof}

Next, let us show that the polygon $Q$ is strictly convex. Given 4 indices
$i_1,i_2,i_3,i_4$, consider two lines in $\bbRP^2$: $(V_{i_1},V_{i_2})$ and
$(V_{i_3},V_{i_4})$; the union of these lines is a conic that does not intersect $P$
anymore. The corresponding hyperplane in~$\bbRP^5$ contains the vertices 
${\cal V}(V_{i_1}),{\cal V}(V_{i_2}),{\cal V}(V_{i_3}),{\cal V}(V_{i_4})$ and intersects $Q$
with multiplicity 4.

Theorem \ref{SixTheorem} now follows from Theorem \ref{BarnerThm} for $d=5$, provided
$Q$ is in general position. If not, then replace
$P$ by a convex polygon $P'$, close to $P$, such that the corresponding polygon $Q'$ is in
general position. Then, as above, $P'$ has at least 6 extremal quintuples of vertices, and
therefore so does $P$. This completes the proof. 

\medskip

To prove Theorems \ref{FourTheorem} and \ref{GhysTheorem}, one replaces the Veronese
map (\ref{VeroneseMap}) by the Veronese map ${\cal V}:\bbRP^2\to\bbRP^3$
$$
{\cal V}:(x:y:z)\mapsto
(x^2+y^2:z^2:yz:zx)
$$
and the Segre map ${\cal S}:\bbRP^1\times\bbRP^1\to\bbRP^3$
$$
{\cal S}:((x_1:y_1),(x_2:y_2))\mapsto
(x_1x_2:x_1y_2:y_1x_2:y_1y_2)
$$
respectively. The proofs of strict convexity of the corresponding polygons $Q$ repeat those
in the smooth case (see \cite{T1}).

\subsection{Conjectures}\label{Discussion}

In conclusion, we formulate three conjectures. Each one is a discrete analog of a theorem
known in the smooth case. We are confident that these conjectures hold true; it would be
interesting to find specifically discrete proofs.

\begin{con}
\label{MobCon}
An embedded non-contractible closed polygon in $\bbRP^2$ has at least 3 flattenings.
\end{con}
\noindent 
In the smooth case this is a celebrated M\"obius theorem (in dimension 2 ``flattening''
means ``inflection'').

\medskip

The notion of flattening of a polygon line extends, in a obvious way, from $\bbRP^d$ to the
sphere $S^d$. 
\begin{con}
\label{SegreCon}
An embedded closed polygon in $S^2$ bisecting the area has at least 4 flattenings.
\end{con}
\noindent 
In the smooth case this was proved by B. Segre \cite{Se} and V. Arnold \cite{A1,A3}.

\medskip

For $k\geq{}d+2$ the next statement is stronger than Theorem \ref{BarnerThm}. 
\begin{con}
\label{BarnerCon}
A strictly convex polygon in $\bbRP^d$ that intersects a hyperplane with multiplicity
$k$ has at least $k$ flattenings.
\end{con}
\noindent 
In the smooth
case this is precisely Barner's result in full generality \cite{Ba}. Conjecture
\ref{BarnerCon} would imply strengthenings of Theorems \ref{FourTheorem},
\ref{SixTheorem} and \ref{GhysTheorem} -- see \cite{T1} for the smooth case.
For example, the following statement would hold.

{\it 
Let $X$ and $Y$ be two $n$-tuples of points in $\bbRP^1$ (see Section \ref{GhysSection}).
If the closed broken line $((x_1,y_1),(x_2,y_2),\ldots,(x_n,y_n))$ in
$\bbRP^1\times\bbRP^1$ intersects the graph of a projective transformation with
multiplicity $k$, then there exists at least $k$ extremal triples of indices.
}

\medskip

\textbf{Acknowledgments:} This work was supported by the Volkswagen-Stiftung
(RiP-program at Oberwolfach). We are grateful to Mathematisches
Forschungsinstitut at Oberwolfach for the creative atmosphere. The second author is also
grateful to Max-Planck-Institut at Bonn for its hospitality. The second author was
supported by an NSF grant.



\begin{thebibliography}{99}


\bibitem{A1} 
V. Arnold, {\it Topological Invariants of Plane Curves and Caustics},
University Lecture Series 5, AMS 1994.

\bibitem{A3} V. Arnold, {\it Topological problems of the theory of wave propagation}, 
Russ. Math. Surv., {\bf 51:1} (1996) 1--47.

\bibitem{Ba} 
M. Barner, {\it Uber die Mindestanzahl stationarer Schmiegebenen bei geschlossenen
strengconvexen Raumcurven}, Abh. Math. Sem. Univ. Hamburg, {\bf 20} (1956) {196--215}.

\bibitem{Bl} 
W. Blaschke, {\it Vorlesungen \"uber Differentielgeometrie}, vol.~2, Springer-Verlag~1923.

\bibitem{C} 
R. Connely,
{\it Rigidity}, in: Handbook of convex geometry, North-Holland, 1993, 223--272.

\bibitem{DO} 
C. Duval \& V. Ovsienko,
{\em 	Schwarzian derivative and Lorentzian wordlines},
Funct. Analysis Appl., to appear.

\bibitem{G} 
E. Ghys, {\it Cercles osculateurs et g\'eom\'etrie lorentzienne}, Talk at the Journ\'ee
Inaugurale du CMI, 1995, Marseille.

\bibitem{GMO} 
 L. Guieu, E. Mourre \& V.Ovsienko,	
{\em 	Theorem on 
Six vertices of a plane curve via the Sturm theory},
Arnold-Gelfand Mathematical Seminars,
Birkh\"auser, 1997,  257--266.

\bibitem{M}
S. Mukhopadhyaya, {\it New methods in the geometry of plane arc},
Bull. Calcutta Math. Soc. {\bf 1} (1909) 31--37.

\bibitem{OT}
V. Ovsienko \& S. Tabachnikov,
{\it	Sturm theory, Ghys theorem on zeroes
of the Schwarzian derivative and flattening of Legendrian curves},
Selecta Mathematica (N. S.),
{\bf 2:2} (1996) 297--307.

\bibitem{S}
V. Sedykh,
{\it Theorem of four support vertices of a polygonal line},
Funct. Analysis Appl., {\bf 30:3} (1996) 216--218.

\bibitem{S1}
V. Sedykh,
{\it Discrete versions of four-vertex theorem},
AMS Transl. Ser.~2, vol.~180, 1997,
197--207.



\bibitem{Se}
B. Segre,
{\it Alcune propriet\`a differenziali in grande delle curve chiuse sghembe},
Rend. Mat. {\bf 1:6} (1968) 237--297.

\bibitem{T1}
S. Tabachnikov,
{\it	On zeroes
of the Schwarzian derivative},
AMS Transl. Ser.~2, vol.~180, 1997, 229--239.

\bibitem{T2}
S. Tabachnikov,
{\it A four vertex theorem for polygons},
Preprint MPIM-Bonn, 1999.

\bibitem{TU}
G. Thorbergsson \& M. Umehara,
{\it A unified approach to the four vertex theorem II}, 
AMS Transl. Ser.~2, vol.~190, 1999, 185--228.

\bibitem{U}
M. Umehara,
{\it A unified approach to the four vertex theorem I}, 
AMS Transl. Ser.~2, vol.~190, 1999, 229--252.

\bibitem{W}
B. Wegner,
{\it On the evolutes of piecewise linear curves in the plane}.
Rad. Hrvat. Acad. Znan. Um. {\bf 467} (1994) 1--16.

\end{thebibliography}
\end{document}